\documentclass[11pt]{article}
\usepackage{amssymb,amsfonts,amsmath,amsthm}
\usepackage{epsfig}
\newtheorem{thm}{Theorem}[section]
\newtheorem{cor}{Corollary}
\newtheorem{lem}[thm]{Lemma}

\makeatletter \@addtoreset{equation}{section}

\title{\bf\LARGE The Bivariate Rogers-Szeg\"{o} Polynomials}
\author{William Y. C. Chen$^{1}$,
Husam L. Saad$^{2}$, and Lisa H. Sun$^{3}$}
\date{$^{1,2,3}$Center for Combinatorics, LPMC,\\
 Nankai
University\\
Tianjin 300071, P.R. China\\
\vskip 0.2 cm $^2$Department of Mathematics, College of Science,\\
Basrah University, Basrah, Iraq\vskip 0.2 cm $^1$chen@nankai.edu.cn,
$^2$hus6274@hotmail.com, $^3$sun@cfc.nankai.edu.cn}

\begin{document}
\maketitle

\noindent \textbf{Abstract.} We present an operator approach to
deriving Mehler's formula and the Rogers formula for the bivariate
Rogers-Szeg\"{o} polynomials $h_n(x,y|q)$. The proof of Mehler's formula
can be considered as a new approach to the
nonsymmetric Poisson kernel formula for the continuous
big $q$-Hermite polynomials $H_n(x;a|q)$ due to Askey, Rahman and Suslov.
Mehler's formula for $h_n(x,y|q)$ involves a
${}_3\phi_2$ sum and the Rogers formula involves a ${}_2\phi_1$
sum. The proofs of these results are based on parameter
augmentation with respect to the $q$-exponential operator and the
homogeneous $q$-shift operator in two variables. By extending
recent results on the Rogers-Szeg\"{o} polynomials $h_n(x|q)$ due to Hou, Lascoux and Mu,
 we
obtain another Rogers-type formula for $h_n(x,y|q)$. Finally, we
give a change of base formula for $H_n(x;a|q)$ which can be used to evaluate
some integrals by using the Askey-Wilson integral.




\section{Introduction}

The Rogers-Szeg\"{o} polynomials $h_n(x|q)$ have been extensively studied since the end of the
nineteenth century. Two classical results for the Rogers-Szeg\"{o} polynomials
are Mehler's formula and the Rogers formula which  respectively correspond to the Poisson
kernel formula and the linearization formula.
In this paper, we extend  Mehler's formula and the Rogers formula to the
bivariate Rogers-Szeg\"{o} polynomials $h_n(x,y|q)$ by using the $q$-exponential
operator as studied in \cite{chenliuII97} and the homogeneous $q$-shift operator
recently introduced by Chen, Fu and Zhang \cite{chenfuzhang}. It should be noted
 that Mehler's formula for $h_n(x,y|q)$ is equivalent to the nonsymmetric Poisson kernel formula
for the continuous big $q$-Hermite polynomials due to Askey, Rahman and
Suslov \cite{askey-rahman-suslov}. So our proof of Mehler's formula for $h_n(x,y|q)$
may be considered as a new approach to the nonsymmetric Poisson kernel for
the continuous big $q$-Hermite polynomials.
As can be seen, the bivariate version of the Rogers-Szeg\"o polynomials is easier to
deal with from the operator point of view.

Let us review some common notation and terminology for basic
hypergeometric series in \cite{gasperrahman90}. Throughout this
paper, we assume that $|q|<1$. The $q$-shifted factorial is
defined by
\begin{equation*}
(a;q)_0=1,\quad (a;q)_\infty=\prod_{k=0}^\infty (1-aq^k),\quad (a;
q)_n = \prod_{k=0}^{n-1}(1-aq^k),\ n\in \mathbb{Z}.
\end{equation*}
The following notation stands for the multiple $q$-shifted
factorials:
\begin{align*}
(a_1,a_2,\ldots, a_m;q)_n =&\  (a_1;q)_n(a_2;q)_n\cdots(a_m;q)_n,\\
(a_1,a_2,\ldots, a_m;q)_\infty =&\
(a_1;q)_\infty(a_2;q)_\infty\cdots(a_m;q)_\infty.
\end{align*}
The $q$-binomial coefficients, or the Gaussian coefficients, are
given by
$$
{n\brack k}=\frac{(q;q)_n}{(q;q)_k(q;q)_{n-k}}.
$$
The basic hypergeometric series ${}_{r+1}\phi_r$ are defined by
$$
\ _{r+1}\phi_r\left(\begin{array}{c}
  a_1, \ldots, a_{r+1}\\
  b_1,\ldots, b_r\\
  \end{array}; q, x \right)=\sum_{n=0}^\infty
  \frac{(a_1,\ldots, a_{r+1};q)_n}{(q,b_1,\ldots,b_r;q)_n}x^n.
$$

We will be mainly concerned with the {\it bivariate
Rogers-Szeg\"{o} polynomials} as given below
$$
h_n(x,y|q)=\sum_{k=0}^n {n\brack k}P_k(x,y),
$$
where $P_n(x,y)=(x-y)(x-qy)\cdots(x-q^{n-1}y)$ are the Cauchy
polynomials with the generating function
\begin{equation}\label{gf-pn}
\sum_{n=0}^\infty P_n(x,y)
\frac{t^n}{(q;q)_n}=\frac{(yt;q)_\infty}{(xt;q)_\infty}, \quad
|xt|<1.
\end{equation}
Note that the Cauchy polynomials $P_{n}(x,y)$ naturally  arise in
the $q$-umbral calculus as studied by Andrews
\cite{andrews71,andrews74}, Goldman and Rota
\cite{Goldman-Rota70}, Goulden and Jackson
\cite{Goulden-Jackson83}, Ihrig and Ismail \cite{Ihrig-Ismail81},
Johnson \cite{Johnson95}, and Roman \cite{Roman82}. The generating
function (\ref{gf-pn})  is also the homogeneous version of the
Cauchy identity, or the $q$-binomial theorem \cite{gasperrahman90}:
\begin{equation}\label{cauchyidentity}
\sum_{k=0}^\infty
\frac{(a;q)_k}{(q;q)_k}z^k=\frac{(az;q)_\infty}{(z;q)_\infty},\quad
|z|<1.
\end{equation}
Putting
$a=0$, \eqref{cauchyidentity} becomes Euler's identity
\cite{gasperrahman90}
\begin{equation}\label{euler-identity}
\sum_{n=0}^\infty \frac{z^k}{(q;q)_k}=\frac{1}{(z;q)_\infty},
\quad |z|<1,
\end{equation}
and its inverse relation takes the form
\begin{equation}\label{inverse-euler}
\sum_{k=0}^\infty \frac{(-1)^kq^{k \choose
2}z^k}{(q;q)_k}=(z;q)_\infty.
\end{equation}

The continuous big $q$-Hermite polynomials
\cite{Koekoek-Swarttouw}  are defined by
$$ H_n(x;a|q)=\sum_{k=0}^n
{n\brack k}(ae^{i\theta};q)_k e^{i(n-2k)\theta},\quad
x=\cos\theta.
$$

We first observe that the bivariate Rogers-Szeg\"{o} polynomials
$h_n(x,y|q)$ introduced by Chen, Fu and Zhang \cite{chenfuzhang}
are equivalent to  the continuous big $q$-Hermite polynomials
owing to the following relation:
\begin{equation}\label{Hxy-hxy}
H_n(x;a|q)=e^{in\theta}h_n(e^{-2i\theta},ae^{-i\theta}|q), \quad
x=\cos\theta.
\end{equation}

The polynomials $h_n(x,y|q)$ have the generating function
\cite{chenfuzhang}:
\begin{equation}\label{gf-hxy}
\sum_{n=0}^{\infty}h_n(x,y|q)\frac{t^n}{(q;q)_n}=
\frac{(yt;q)_\infty}{(t,xt;q)_\infty}, \quad |t|<1, |xt|<1,
\end{equation}
which is equivalent to the generating function for the big
continuous $q$-Hermite polynomials, see, for example,
Koekoek-Swarttouw \cite{Koekoek-Swarttouw}.  Notice that the
classical {\it Rogers-Szeg\"{o} polynomials}
$$
h_n(x|q)=\sum_{k=0}^n {n\brack k}x^k,
$$
are a special case of $h_n(x,y|q)$ when $y$ is set to zero, and in
this case (\ref{gf-hxy}) reduces to
\begin{equation}\label{gf-hx}
\sum_{n=0}^\infty
h_n(x|q)\frac{t^n}{(q;q)_n}=\frac{1}{(t,xt;q)_\infty},\quad |t|<1.
\end{equation}
The Rogers-Szeg\"{o} polynomials play an important role in the
theory of orthogonal polynomials, particularly in the study of the
Askey-Wilson polynomials, see
\cite{Salam-Ismai88,askeyismail83,Atkishiyev-Nagiyev94,
bressoud80,ismail-stanton88, i-s-v87, Rogers1893,stanton2000}.
They are closely related to the $q$-Hermite polynomials
$$
H_n(x|q)=\sum_{k=0}^n {n\brack k}e^{i(n-2k)\theta},\quad x=\cos
\theta.
$$
In fact,  the following relations hold
\begin{equation}\label{Hx-hx}
H_n(x|q)=H_n(x;0|q)=e^{in\theta}h_n(e^{-2i\theta}|q), \quad
x=\cos\theta.
\end{equation}
The continuous big $q$-Hermite polynomials $H_n(x;a|q)$ are connected with
 the $q$-Hermite polynomials
$H_n(x|q)$ via the following relation \cite{AtaAta97, FLV95}:
\begin{equation}\label{Hxa-Hx}
H_n(x;a|q)=\sum_{k=0}^n{n\brack k}(-1)^kq^{{k\choose
2}}a^kH_{n-k}(x|q),
\end{equation}
and the inverse expansion of \eqref{Hxa-Hx} becomes
\begin{equation}\label{Hx-Hxa}
H_n(x|q)=\sum_{k=0}^n {n\brack k}a^kH_{n-k}(x;a|q).
\end{equation}

This paper is motivated by the natural question of extending
Mehler's formula to $h_n(x,y|q)$, where
Mehler's formula for the Rogers-Szeg\"{o} polynomials reads
\begin{equation}\label{mehler-hx}
\sum_{n=0}^{\infty}h_{n}(x|q)h_{n}(y|q)
\frac{t^{n}}{(q;q)_{n}}=\frac{(xyt^2;q)_{\infty}}{(t,xt,yt,xyt;q)_{\infty}}.
\end{equation}
The formula (\ref{mehler-hx}) has been extensively studied, see
\cite{chenliuII97, ismail-stanton88, Karande-Thakare74,
Rogers1893, stanton2000, Zhang-Liu06}. Based on the recurrence
relation for $H_n(x|q)$, Bressoud \cite{bressoud80} gave a proof
of the equivalent formula, or the Poisson kernel formula, for the $q$-Hermite
polynomials $H_n(x|q)$. Ismail, Stanton and Viennot \cite{i-s-v87} found
a combinatorial proof of the Poisson kernel formula for $H_n(x|q)$ by using the vector
space interpretation of the $q$-binomial coefficients. Askey, Rahman and
Suslov \cite{askey-rahman-suslov} derived the nonsymmetric Poisson kernel
formula for the continuous big $q$-Hermite polynomials
$H_n(x;a|q)$:
\begin{align}\label{nonsympoiker}
\sum_{n=0}^\infty H_n(x;a|q)H_n(y;b|q)\frac{t^n}{(q;q)_n}=&\
\frac{(ate^{i\beta},be^{-i\beta}, t^2;q)_\infty}{(te^{i(\theta+
\beta)},te^{i(\theta-\beta)},te^{-i(\theta+\beta)},
te^{-i(\theta-\beta)};q)_\infty}\nonumber\\
&\times\ _3\phi_2\left(\begin{array}{c}
  t e^{i(\theta+\beta)}, t e^{-i(\theta-\beta)}, at/b\\
  ate^{i\beta}, t^2\\
  \end{array}; q, be^{-i\beta}\right),
\end{align}
where $x=\cos\theta,\ y=\cos \beta$. The above formula can be viewed as
Mehler's formula for $H_n(x;a|q)$. Moreover, it can be restated in terms
of $h_n(x,y|q)$  (Theorem \ref{thm-mehler-hxy}). The first result of this paper
is an operator approach to Mehler's formula for $h_n(x,y|q)$. We will
present a simple proof by using the exponential operators involving the
classical $q$-differential operator and a bivariate $q$-differential operator
introduced by Chen, Fu and Zhang \cite{chenfuzhang}.

The second result of this paper is the Rogers formula for
$h_n(x,y|q)$. The Rogers formula \cite{chenliuII97,
Rogers1893,Rogers1893II} for the classical Rogers-Szeg\"{o}
polynomials $h_n(x|q)$ reads:
\begin{align}\label{rogers-hx}
\sum_{n=0}^{\infty}\sum_{m=0}^{\infty}h_{n+m}(x|q)
\frac{t^{n}}{(q;q)_{n}} \frac{s^{m}}{(q;q)_{m}}
=(xst;q)_{\infty}\sum_{n=0}^{\infty}
\sum_{m=0}^{\infty}h_{n}(x|q)h_{m}(x|q) \frac{t^{n}}{(q;q)_{n}}
\frac{s^{m}}{(q;q)_{m}}.
\end{align}
One of the most important applications of the Rogers formula is to deduce
the following linearization formula for $h_n(x|q)$ (cf.  Bressound \cite{bressoud80},
Ismail-Stanton \cite{ismail-stanton88} and Rogers
\cite{rogers1895})
\begin{equation}\label{linearhx}
h_{n}(x|q)h_{m}(x|q)=\sum_{k=0}^{\min\{n,m\}}{n \brack k}{m \brack
k}(q;q)_{k}x^{k}h_{n+m-2k}(x|q).
\end{equation}

Based on a recent approach of Hou, Lascoux and Mu to the Rogers-Szeg\"o polynomials,
we derive a second Rogers-type formula for $h_n(x,y|q)$ which leads to a simpler
linearization formula compared with the first one we have obtained.

We conclude this paper with a change of base formula for  $H_n(x;a|q)$.
This formula along with other identities can be used to compute some integrals
with the aid of the Askey-Wilson integral.

\section{Mehler's Formula for $h_n(x,y|q)$}

In this section, we aim to present an operator approach to Mehler's formula
for the bivariate Rogers-Szeg\"o polynomials $h_n(x,y|q)$.

\begin{thm}[Mehler's Formula for
$h_n(x,y|q)$] We have \label{thm-mehler-hxy}
\begin{equation}\label{mehler-hxy}
\sum_{n=0}^\infty h_n(x,y|q)h_n(u,v|q)\frac{t^n}{(q;q)_n}=
\frac{(yt,vxt;q)_\infty}{(t,xt,uxt;q)_\infty}\
{}_3\phi_2\left(\begin{array}{c}
  y, xt, v/u\\
  yt, vxt\\
  \end{array}; q, ut\right),
\end{equation}
provided that $|t|, |xt|, |ut|, |uxt|<1$.
\end{thm}
Obviously, Mehler's formula (\ref{mehler-hx}) for $h_n(x|q)$ can
be deduced from the above theorem by  setting $y=0, v=0$ and
$u=y$. We note that it is not difficult to reformulate
\eqref{mehler-hxy} as the nonsymmetric Poisson kernel formula
\eqref{nonsympoiker} for $H_n(x;a|q)$. To this end, we first make
the variable substitutions $x\rightarrow e^{-2i\theta},
y\rightarrow
    ae^{-i\theta}, u\rightarrow e^{-2i\beta}, v\rightarrow
    be^{-i\beta}$ so that
we may use the relation (\ref{Hxy-hxy}) to transform $h_n(x,y|q)$
and $h_n(u,v|q)$ into $H_n(x;a|q)$
    and $H_n(y;b|q)$. Then the formula \eqref{nonsympoiker} follows
from  the $_3\phi_2$ transformation \cite[Appendix III.
9]{gasperrahman90}:
\begin{align}\label{3phi2trans}
{}_3\phi_2\left(\begin{array}{c}
  a, b, c\\
  d, e\\
  \end{array}; q, \frac{de}{abc}\right)=\frac{(e/a,de/bc;q)_\infty}{(e,de/abc;q)_\infty}
  \ {}_3\phi_2\left(\begin{array}{c}
  a, d/b, d/c\\
  d, de/bc\\
  \end{array}; q, \frac{e}{a}\right).
\end{align}

Our operator approach to Theorem \ref{mehler-hxy} involves
two identities (Lemmas
\ref{lem22} and \ref{lem23}) in connection with the
$q$-exponential operator and the homogeneous $q$-shift operator.
The $q$-differential operator, or the $q$-derivative, acting on
the variable $a$, is defined by
$$
D_qf(a)=\frac{f(a)-f(aq)}{a},
$$
and the $q$-exponential operator is given by
$$
T(bD_q)=\sum_{n=0}^\infty \frac{(bD_q)^n}{(q;q)_n}.
$$
Evidently,
\begin{equation}\label{tdq-xn}
T(D_{q})\{x^{n}\}=h_{n}(x|q).
\end{equation}

\begin{lem} \label{lem22}  We have
\begin{equation}\label{eq25}
T(bD_q)\left\{\frac{(av;q)_\infty}{(as,at;q)_\infty}\right\}=
\frac{(bv;q)_\infty}{(as,bs,b\,t;q)_\infty}\
_2\phi_1\left(\begin{array}{c}
  v/t, bs\\
  bv\\
  \end{array}; q, at\right),
\end{equation}
provided that $|bs|, |b\,t|<1$.
\end{lem}

From the Leibniz rule for $D_q$ (see \cite{Roman85})
\[
D_q^n\{f(a)g(a)\}=\sum_{k=0}^n q^{k(k-n)}{n \brack k}D_q^k\{f(a)\}D_q^{n-k}\{g(q^ka)\},
\]
\eqref{eq25} can be verified by straightforward computation.
Here we also note that a more general relation has been established by
Zhang and Wang \cite{zhangwang05}:
\begin{align}\label{tdq3phi2}
T(bD_q)\left\{\frac{(av;q)_\infty}{(as,at,aw;q)_\infty}\right\}
=&\ (av,bv;q)_\infty
\frac{(abstw/v;q)_\infty}{(as,at,aw,bs,b\,t,bw;q)_\infty}\nonumber\\
&\ \times\, {}_3\phi_2\left(\begin{array}{c}
  v/s, v/t, v/w\\
  av, bv\\
  \end{array}; q, abstw/v\right),
\end{align}
where $|bs|, |b\,t|, |bw|, |abstw/v|<1$. Setting $w=0$ in
(\ref{tdq3phi2}), by virtue of Jackson's transformation \cite[Appendix III.
4]{gasperrahman90} and Heine's transformation \cite[Appendix III.
1]{gasperrahman90}, \eqref{eq25} becomes a consequence of (\ref{tdq3phi2}).

In \cite{chenfuzhang}, Chen, Fu and Zhang introduced the homogeneous
$q$-difference operator
$$
D_{xy}f(x,y)=\frac{f(x,q^{-1}y)-f(qx,y)}{x-q^{-1}y}
$$
and the homogeneous $q$-shift operator
\[
\mathbb{E}(D_{xy})=\sum_{k=0}^\infty \; {D_{xy}^k\over (q;q)_k}.
\]
The following basic facts have been observed in
\cite{chenfuzhang}:
\begin{align}
& D_{xy}\{P_n(x,y)\}=(1-q^n)P_{n-1}(x,y),\nonumber\\[6pt]
&\mathbb{E}(D_{xy})\{P_n(x,y)\}=h_n(x,y|q).\label{Edxy-pn}
\end{align}

\begin{lem}\label{lem23}
We have
$$
\mathbb{E}(D_{xy})\left\{\frac{(yt;q)_\infty}{(xt;q)_\infty}
\frac{P_n(x,y)}{(yt;q)_n}\right\}=\frac{(yt;q)_\infty}{(t,xt;q)_\infty}
\sum_{k=0}^n {n\brack k}\frac{(y,xt;q)_k}{(yt;q)_k}x^{n-k},
$$
provided that $|t|, |xt|<1$.
\end{lem}

\begin{proof}
Let us compute the following sum in two ways:
\begin{equation}\label{hnxyz}
\sum_{n=0}^{\infty}h_{n}(x,y|q)h_{n}(z|q)\frac{t^{n}}{(q;q)_{n}}.
\end{equation}
We may either express $h_n(z|q)$ as
$T(D_q)\{z^n\}$ by \eqref{tdq-xn} or express $h_n(x,y|q)$ as
$\mathbb{E}(D_{xy})\{P_n(x,y)\}$ by \eqref{Edxy-pn}.
Invoking $h_n(z|q)=T(D_q) \{z^n\}$, the sum \eqref{hnxyz} equals
\begin{align*}
\lefteqn{\sum_{n=0}^{\infty}h_{n}(x,y|q)T(D_q)\left\{z^n\right\}
\frac{t^{n}}{(q;q)_{n}}}\\
&\ =T(D_q)\left\{\sum_{n=0}^{\infty}h_{n}(x,y|q)
\frac{(zt)^{n}}{(q;q)_{n}}\right\}\quad (|zt|<1,\ |xzt|<1)\\
&\ =T(D_q)\left\{\frac{(yzt;q)_\infty}{(xzt,zt;q)_\infty}\right\}(|t|<1,\ |xt|<1).
\end{align*}
According to Lemma \ref{lem22}, \eqref{hnxyz} can be expressed in the following form
\begin{equation}
\frac{(yt;q)_{\infty}}{(xzt,xt,t;q)_{\infty}}\
_2\phi_1\left(\begin{array}{c}
  y, xt\\
  yt\\
  \end{array}; q, zt\right).\label{rhs}
\end{equation}
On the other hand, \eqref{hnxyz} also equals
\begin{align*}
\lefteqn{\sum_{n=0}^{\infty}\mathbb{E}(D_{xy})\left\{P_{n}(x,y)\right\}h_{n}(z|q)
\frac{t^{n}}{(q;q)_{n}}}\nonumber\\
&\ =\mathbb{E}(D_{xy})\left\{\sum_{n=0}^{\infty}P_n(x,y)h_{n}(z|q)
\frac{t^{n}}{(q;q)_{n}}\right\}\nonumber\\
&\ =\mathbb{E}(D_{xy})\left\{\sum_{n=0}^{\infty}P_n(x,y)\sum_{k=0}^n
{n\brack k}z^k\frac{t^{n}}{(q;q)_{n}} \right\}\nonumber\\
&\ =\mathbb{E}(D_{xy})\left\{\sum_{k=0}^{\infty}\left(\sum_{n=0}^{\infty}P_n(x,q^{k}y)
\frac{t^{n}}{(q;q)_{n}}\right)P_k(x,y) \frac{(zt)^{k}}{(q;q)_{k}}
\right\}\nonumber\\
&\ =\sum_{k=0}^{\infty}\frac{(zt)^k}{(q;q)_k}\mathbb{E}(D_{xy})
\left\{\frac{(yt;q)_\infty}{(xt;q)_\infty}
\frac{P_k(x,y)}{(yt;q)_k}\right\},
\end{align*}
where $|xt|<1$. Now we see that
$$
\sum_{k=0}^{\infty}\frac{(zt)^k}{(q;q)_k}\mathbb{E}(D_{xy})
\left\{\frac{(yt;q)_\infty}{(xt;q)_\infty}
\frac{P_k(x,y)}{(yt;q)_k}\right\}=\frac{(yt;q)_{\infty}}{(xzt,xt,t;q)_{\infty}}\
_2\phi_1\left(\begin{array}{c}
  y, xt\\
  yt\\
  \end{array}; q, zt\right).
$$
Employing Euler's identity
\eqref{euler-identity} for $1/(xzt;q)_\infty$ and expanding the ${}_2\phi_1$ summation
on the right hand side of the
above identity, we obtain
\begin{align*}
\sum_{k=0}^{\infty}\frac{(zt)^k}{(q;q)_k}\mathbb{E}(D_{xy})
\left\{\frac{(yt;q)_\infty}{(xt;q)_\infty}
\frac{P_k(x,y)}{(yt;q)_k}\right\}=\frac{(yt;q)_{\infty}}
{(t,xt;q)_{\infty}}\sum_{n=0}^{\infty}\sum_{k=0}^{\infty}
  \frac{(y,xt;q)_n}{(q,yt;q)_n}\frac{z^{n+k}t^{n+k}x^{k}}{(q;q)_{k}}.
\end{align*}
Equating the coefficients of $z^n$,  the desired identity follows.
\end{proof}

We are now ready to present the proof of  Theorem
\ref{thm-mehler-hxy}.

\begin{proof} From (\ref{Edxy-pn}) it
follows that
\begin{align*}
\lefteqn{\sum_{n=0}^\infty h_n(x,y|q)h_n(u,v|q)\frac{t^n}{(q;q)_n}}\\
&\ = \mathbb{E}(D_{xy})\left\{\sum_{n=0}^\infty
P_n(x,y)h_n(u,v|q)\frac{t^n}{(q;q)_n}\right\}\\
&\ = \mathbb{E}(D_{xy})\left\{\sum_{n=0}^\infty
P_n(x,y)\frac{t^n}{(q;q)_n}\sum_{k=0}^n{n\brack k}P_k(u,v)\right\}\\
&\ = \mathbb{E}(D_{xy})\left\{\sum_{k=0}^\infty
P_k(u,v)P_k(x,y)\frac{t^k}{(q;q)_k}\left(\sum_{n=0}^\infty P_n(x,
q^ky)\frac{t^n}{(q;q)_n}\right)\right\}\quad (|xt|<1)\\
&\ = \mathbb{E}(D_{xy})\left\{\sum_{k=0}^\infty
P_k(u,v)P_k(x,y)\frac{t^k}{(q;q)_k}\frac{(q^kyt;q)_\infty}{(xt;q)_\infty}
\right\}\\
&\ =\sum_{k=0}^\infty
P_k(u,v)\frac{t^k}{(q;q)_k}\mathbb{E}(D_{xy})\left\{\frac{(yt;q)_\infty}{(xt;q)_\infty}
\frac{P_k(x,y)}{(yt;q)_k}
\right\}\quad (|t|, |xt|<1).\\
\end{align*}
In view of Lemma \ref{lem23}, the above summation equals
\begin{align*}
\frac{(yt;q)_\infty}{(t,xt;q)_\infty} \sum_{k=0}^\infty
P_k(u,v)\frac{t^k}{(q;q)_k} \sum_{j=0}^k {k\brack
j}\frac{(y,xt;q)_j}{(yt;q)_j}x^{k-j}.
\end{align*}
Exchanging the order of summations, we get
\begin{align*}
 \lefteqn{\frac{(yt;q)_\infty}{(t,xt;q)_\infty}
\sum_{j=0}^\infty P_j(u,v)\frac{(y,xt;q)_j}{(q,yt;q)_j}t^j
\sum_{k=0}^\infty \frac{(xt)^kP_k(u,q^jv)}{(q;q)_k}\quad
(|uxt|<1)}\\
&\ =\frac{(yt,vxt;q)_\infty}{(t,xt,uxt;q)_\infty}
\sum_{j=0}^\infty
P_j(u,v)\frac{(y,xt;q)_j}{(q,yt,vxt;q)_j}t^j\\
&\ = \frac{(yt,vxt;q)_\infty}{(t,xt,uxt;q)_\infty}\
_3\phi_2\left(\begin{array}{c}
  y, xt, v/u\\
  yt, vxt\\
  \end{array}; q, ut\right)\quad (|ut|<1).
  \end{align*}
This completes the proof.
\end{proof}

\section{The Rogers Formula for $h_n(x,y|q)$}

In this section, we obtain the Rogers formula for the bivariate Rogers-Szeg\"{o} polynomials
$h_n(x,y|q)$ using the operator $\mathbb{E}(D_{xy})$ and the
technique of parameter augmentation \cite{chenfuzhang,
chenliuII97}. This Rogers formula implies a linearization formula for $h_n(x,y|q)$.
We also get another Rogers-type formula for $h_n(x,y|q)$
which leads to a simpler linearization formula.

\begin{thm}[\bf{The Rogers Formula for $h_n(x,y|q)$}]\label{rogersformula}
We have
\begin{align}\label{rogersforhxy}
&\ \sum_{n=0}^\infty\sum_{m=0}^\infty
h_{n+m}(x,y|q)\frac{t^n}{(q;q)_n}\frac{s^m}{(q;q)_m} =
\frac{(ys;q)_\infty}{(s,xs,xt;q)_\infty}\
_2\phi_1\left(\begin{array}{c}
  y, xs\\
  ys\\
  \end{array}; q, t\right),
\end{align}
provided that $|t|, |s|, |xt|, |xs|<1$.
\end{thm}

\begin{proof}  By (\ref{Edxy-pn}), we have
\begin{align*}
\lefteqn{\sum_{n=0}^\infty\sum_{m=0}^\infty
h_{n+m}(x,y|q)\frac{t^n}{(q;q)_n}\frac{s^m}{(q;q)_m}}\\
&\ = \mathbb{E}(D_{xy})\left\{\sum_{n=0}^\infty\sum_{m=0}^\infty
P_{n+m}(x,y)\frac{t^n}{(q;q)_n}\frac{s^m}{(q;q)_m}\right\}\\
&\ = \mathbb{E}(D_{xy})\left\{\sum_{n=0}^\infty
P_n(x,y)\frac{t^n}{(q;q)_n}\left(\sum_{m=0}^\infty
P_m(x,q^ny)\frac{s^m}{(q;q)_m}\right)\right\}\quad (|xs|<1)\\
&\ = \mathbb{E}(D_{xy})\left\{\sum_{n=0}^\infty
P_n(x,y)\frac{t^n}{(q;q)_n}\frac{(q^nys;q)_\infty}{(xs;q)_\infty}\right\}\\
&\ = \sum_{n=0}^\infty\frac{t^n}{(q;q)_n}\mathbb{E}(D_{xy})
\left\{\frac{(ys;q)_\infty P_n(x,y)}{(xs;q)_\infty
(ys;q)_n}\right\}\quad (|s|<1,\ |xs|<1).
\end{align*}
Empolying Lemma \ref{lem23}, we find
\begin{align*}
\lefteqn{
\frac{(ys;q)_\infty}{(s,xs;q)_\infty}\sum_{n=0}^\infty\frac{t^n}{(q;q)_n}
 \sum_{k=0}^n {n\brack
k}\frac{(y,xs;q)_k}{(ys;q)_k}x^{n-k}}\\
&\
=\frac{(ys;q)_\infty}{(s,xs;q)_\infty}\sum_{k=0}^\infty\frac{(y,xs;q)_k}{(q,ys;q)_k}t^{k}
\sum_{n=0}^\infty\frac{(xt)^n}{(q;q)_n}\quad (|xt|<1)\\
&\ = \frac{(ys;q)_\infty}{(s,xs,xt;q)_\infty} \sum_{k=0}^\infty
\frac{(y,xs;q)_k}{(q,ys;q)_k}t^{k}\\
&\ = \frac{(ys;q)_\infty}{(s,xs,xt;q)_\infty}\
_2\phi_1\left(\begin{array}{c}
  y, xs\\
  ys\\
  \end{array}; q, t\right)\quad (|t|<1),
\end{align*}
as desired.
\end{proof}

Clearly, the Rogers formula (\ref{rogers-hx}) for $h_n(x|q)$ is a
special case of (\ref{rogersforhxy}) when $y=0$.  From the above
theorem and \eqref{Hxy-hxy}, we get the equivalent formula for
$H_n(x;a|q)$:
\begin{align*} \sum_{n=0}^\infty\sum_{m=0}^\infty
H_{n+m}(x;a|q)\frac{t^n}{(q;q)_n}\frac{s^m}{(q;q)_m} =&\
\frac{(as;q)_\infty}{(se^{i\theta},se^{-i\theta},te^{-i\theta};q)_\infty}\,{}_2\phi_1
\left(\begin{array}{c}
  ae^{-i\theta}, se^{-i\theta}\\
  as\\
    \end{array}; q, te^{i\theta}\right),
\end{align*}
where $x=\cos \theta$ and $|te^{i\theta}|, |se^{i\theta}|,
|te^{-i\theta}|, |se^{-i\theta}|<1$.

As in the classical case, the Rogers formula can be used to derive
linearization formula. For the bivariate case, we obtain the
linearization formula for  $h_n(x,y|q)$ as a double summation
identity.

\begin{cor}We have
\begin{align*}
\lefteqn{\sum_{k=0}^n\sum_{l=0}^m {n\brack k}{m\brack
l}(y;q)_k(y/x;q)_lx^lh_{n+m-k-l}(x,y|q)}\\
&\quad=\sum_{k=0}^n\sum_{l=0}^m {n\brack k}{m\brack
l}(y;q)_k(y/x;q)_l (xq^k)^lh_{n-k}(x,y|q)h_{m-l}(x,y|q).
\end{align*}
\end{cor}

\begin{proof}
We rewrite Theorem~3.1 in the following form
\begin{align*}
\lefteqn{\frac{(ys;q)_\infty(yt;q)_\infty}{(xs;q)_\infty(t;q)_\infty}\sum_{n=0}^\infty\sum_{m=0}^\infty
h_{n+m}(x,y|q)\frac{t^n}{(q;q)_n}\frac{s^m}{(q;q)_m}}\\
&\quad =\sum_{k=0}^\infty
\frac{(y;q)_k(ysq^k;q)_\infty}{(q;q)_k(xsq^k;q)_\infty}t^k
\sum_{n=0}^\infty\sum_{m=0}^\infty
h_n(x,y|q)h_m(x,y|q)\frac{t^n}{(q;q)_n}\frac{s^m}{(q;q)_m}.
\end{align*}
Expanding $(ys;q)_\infty/(xs;q)_\infty,\
(yt;q)_\infty/(t;q)_\infty,\ (ysq^k;q)_\infty/(xsq^k;q)_\infty$ by
the Cauchy identity \eqref{cauchyidentity}, and equating the
coefficients of $t^{n}s^{m}$, the required formula is justified.
\end{proof}

We note that Hou, Lascoux and Mu \cite{hou-lascoux-mu}
represented the Rogers-Szeg\"{o} polynomials $h_n(x|q)$ as a
special case of the complete symmetric functions. By  computing the Hankel forms, they
obtained the  Askey-Ismail formula (see \cite{askeyismail83, chenliuII97}):
\begin{equation}\label{askey-ismail}
h_{m+n}(x|q)=\sum_{k=0}^{\min\{m,n\}}{n \brack k}{m \brack
k}(q;q)_{k}q^{k\choose 2}(-x)^{k}h_{n-k}(x|q)h_{m-k}(x|q),
\end{equation}
which can be regarded as the inverse relation of the linearization formula
\eqref{linearhx} for the Rogers-Szeg\"o polynomials.
Applying the technique of  Hou, Lascoux and Mu  to the bivariate Rogers-Szeg\"{o}
polynomials $h_n(x,y|q)$, the following relation can be verified:
\begin{equation}\label{rogerszegoII}
\sum_{k=0}^{\min\{m,n\}} {n\brack k}{m\brack
k}(-1)^k(q;q)_kq^{k\choose
2}\big(x^kh_{n-k}(x,y|q)h_{m-k}(x,y|q)-y^kh_{n+m-k}(x,y|q)\big)=0.
\end{equation}
The details of the proof are omitted.
Multiplying the above equation
by
$$
\frac{t^n}{(q;q)_n}\frac{s^m}{(q;q)_m}
$$ and summing over $n$ and
$m$, we get another Rogers-type formula.

\begin{thm} We have
\begin{align}
\lefteqn{\sum_{k=0}^\infty\frac{(-1)^ky^kq^{{k\choose
2}}}{(q;q)_k}\sum_{n=k}^\infty\sum_{m=k}^\infty
h_{n+m-k}(x,y|q)\frac{t^n}{(q;q)_{n-k}}
\frac{s^m}{(q;q)_{m-k}}}\nonumber\\
&\quad =(xst;q)_\infty \sum_{n=0}^\infty\sum_{m=0}^\infty
h_n(x,y|q)h_m(x,y|q)\frac{t^n}{(q;q)_n} \frac{s^m}{(q;q)_m}.
\end{align}
\end{thm}

Clearly, the classical Rogers formula \eqref{rogers-hx} is a
special case when $y=0$.
By equating the coefficients of $t^ns^m$ in the above theorem, we
can derive a simpler linearization formula for $h_n(x,y|q)$.

\begin{cor} For $n, m\geq 0$, we have
\begin{align}\label{linearhxy}
\lefteqn{h_n(x,y|q)h_m(x,y|q)\nonumber}\\
&\quad =\sum_{l=0}^{\min\{m,n\}}\sum_{k=0}^{\min\{m,n\}}{m\brack
l}{n\brack l}{m-l \brack k}{n-l \brack
k}(q;q)_k(q;q)_l(-1)^kx^ly^kq^{{k\choose 2}}h_{n+m-2l-k}(x,y|q).
\end{align}
\end{cor}

The following special case of Theorem~\ref{rogersformula}
for $y=0$ will be useful to verify the relation  between
$h_n(x|q)$ and $h_n(x,y|q)$.

\begin{cor} For $n, m\geq 0$, we have
\begin{align}\label{eq32}
\lefteqn{\sum_{k=0}^{\min\{n,m\}}{n \brack k}{m \brack
k}(q;q)_{k}x^{k}h_{n+m-2k}(x|q)}\nonumber\\
&\quad\quad =\left(\sum_{k=0}^{n}{n \brack k}y^{k}\
h_{n-k}(x,y|q)\right)\left(\sum_{j=0}^{m}{m \brack j}y^{j}\
h_{m-j}(x,y|q)\right).
\end{align}
\end{cor}

\begin{proof} Setting $y=0$ in Theorem \ref{rogersforhxy}, from
 the Cauchy identity (\ref{cauchyidentity})
and  (\ref{gf-hxy})  it follows that
\begin{align}\label{partrogers}
\lefteqn{\displaystyle\sum_{n=0}^{\infty}
\sum_{m=0}^{\infty}h_{n+m}(x|q)
\frac{t^{n}}{(q;q)_{n}}\frac{s^{m}}{(q;q)_{m}} \nonumber}\\
&\ =\frac{1}{(s,xs,xt;q)_\infty}\sum_{k=0}^\infty\frac{(xs;q)_k}{(q;q)_k}t^k\quad (|t|<1)\nonumber\\
&\ =\displaystyle
 \frac{(xst;q)_{\infty}}{(ys,yt;q)_\infty}
 \frac{(yt;q)_{\infty}}{(t,xt;q)_{\infty}} \frac{(ys;q)_{\infty}}
 {(s,xs;q)_{\infty}}
\nonumber\quad (|t|, |s|, |xt|, |xs|<1)\\
&\ =\displaystyle
\frac{(xst;q)_{\infty}}{(ys,yt;q)_\infty}\sum_{n=0}^{\infty}
\sum_{m=0}^{\infty}h_{n}(x,y|q)h_{m}(x,y|q)
\frac{t^{n}}{(q;q)_{n}} \frac{s^{m}}{(q;q)_{m}},
\end{align}
which can be rewritten as
\begin{align*}
\lefteqn{\frac{1}{(xst;q)_{\infty}}\sum_{n=0}^{\infty}\sum_{m=0}^{\infty}h_{n+m}(x|q)\
\frac{t^{n}}{(q;q)_{n}} \frac{s^{m}}{(q;q)_{m}}}\\&\quad =
\frac{1}{(yt,ys;q)_{\infty}}\sum_{n=0}^{\infty}
\sum_{m=0}^{\infty}h_{n}(x,y|q)h_{m}(x,y|q)
\frac{t^{n}}{(q;q)_{n}} \frac{s^{m}}{(q;q)_{m}},
\end{align*}
where $|t|, |s|, |xt|, |xs|<1$.

Assuming that $|xst|, |yt|, |ys|<1$,  we can expand
$1/(xst;q)_{\infty},\ 1/(yt;q)_{\infty}$ and $1/(ys;q)_{\infty}$
by Euler's identity \eqref{euler-identity}. Equating
coefficients of $t^{n}s^{m}$ gives (\ref{eq32}). Since $|t|,
|s|, |xt|,$ $|xs|<1$ and $|xst|, |yt|, |ys|<1$,  we see that $|x|$
and $|y|$ must be finite. This completes the proof.
\end{proof}

When $y=0$, both \eqref{linearhxy} and (\ref{eq32}) reduce to the
well-known linearization formula \eqref{linearhx}.
Setting $m=0$ in (\ref{eq32}), we are led to the following
relation between $h_n(x|q)$ and $h_n(x,y|q)$
\begin{equation}\label{eq34}
h_{n}(x|q)=\sum_{k=0}^{n}{n \brack k}y^{k}\ h_{n-k}(x,y|q),
\end{equation}
which is a special case of a relation of Askey-Wilson
\cite[(6.4)]{askey-wilson}. The inverse relation of \eqref{eq34}
is as follows
\begin{equation}\label{eq35}
h_n(x,y|q)=\sum_{k=0}^{n}{n \brack k}(-1)^kq^{k\choose 2}y^k
h_{n-k}(x|q).
\end{equation}
Note that  (\ref{eq34}) and (\ref{eq35}) are equivalent to the
relations (\ref{Hx-Hxa}) and (\ref{Hxa-Hx}) between $H_n(x|q)$ and
$H_n(x;a|q)$.

In fact, we can go one step further from \eqref{partrogers}. Reformulating
\eqref{partrogers} by multiplying $(ys,yt;q)_\infty$ on both
sides and expanding $(ys;q)_{\infty}, (yt;q)_{\infty}$ and
$(xst;q)_{\infty}$ using Euler's formula \eqref{inverse-euler}, we get
\begin{align*}
\lefteqn{\displaystyle\sum_{n=0}^{\infty}\sum_{m=0}^{\infty}\sum_{j=0}^{\infty}\sum_{k=0}^{\infty}
\frac{q^{{j \choose 2}+{k \choose 2}}
(-y)^{j+k}}{(q;q)_{j}(q;q)_{k}}h_{n+m}(x|q)
 \frac{t^{n+j}}{(q,q)_{n}} \frac{s^{m+k}}{(q,q)_{m}}}
\\
&\quad =\sum_{n=0}^{\infty}\sum_{m=0}^{\infty}\sum_{k=0}^{\infty}
\frac{q^{k \choose
2}(-x)^{k}}{(q;q)_{k}}h_{n}(x,y|q)h_{m}(x,y|q)\frac{t^{n+k}}{(q,q)_{n}}
\frac{s^{m+k}}{(q,q)_{m}}.
\end{align*}
Comparing the coefficients of $t^{n}s^{m}$,
we reach the following identity
\begin{align}
\lefteqn{\displaystyle\sum_{j=0}^{n}\sum_{k=0}^{m}{n \brack j}{m
\brack k}q^{{j \choose 2}+{k \choose 2}}(-y)^{j+k}
h_{n+m-(j+k)}(x|q)\nonumber}\\
&\quad =\sum_{k=0}^{\min\{n,m\}}{n \brack k}{m \brack
k}(q;q)_kq^{k \choose 2}(-x)^{k}\ h_{n-k}(x,y|q)h_{m-k}(x,y|q).
\end{align}
Setting $y=0$ in the above identity, we are led to the Askey-Ismail
formula \eqref{askey-ismail}.

\section{A Change of Base Formula for $H_n(x;a|q)$}

In this section, we give an extension of the $q$-Hermite change of
base formula to the continuous big $q$-Hermite polynomials. The
corresponding statement for $h_n(x,y|q)$ are omitted because we
find that it is more convenient to work with $H_n(x;a|q)$ for this
purpose. This formula can be used to evaluate certain integrals.

In \cite[p.\,7]{ismail-stanton03}, Ismail and Stanton gave
the following $q$-Hermite change of base formula for $H_n(x|p)$:
\begin{equation}\label{changebaseH}
H_n(x|p)=\sum_{j=0}^{n/2} c_{n, n-2j}(p,q) H_{n-2j}(x|q),
\end{equation}
where $x=\cos\theta$ and
\begin{align*}
c_{n,n-2k}(p,q)=\sum_{j=0}^k (-1)^jp^{k-j}q^{j+1\choose 2}{n-2k+j \brack j}_q
 \left({n\brack k-j}_p-p^{n-2k+2j+1}{n \brack k-j-1}_p\right).
\end{align*}

From \eqref{Hxa-Hx}, \eqref{Hx-Hxa} along with the above relation, we obtain
{\it a change of base formula for $H_n(x;a|q)$}:
\begin{equation}\label{changebase}
H_n(x;a|p)=\sum_{j=0}^n d_{n,j,l,m}(p,q) H_{n-j-2l-m}(x;a|q),
\end{equation}
where $x=\cos\theta$ and
$$
d_{n,j,l,m}(p,q)={n \brack j}_p(-1)^jp^{j \choose 2}a^j
\sum_{l=0}^{(n-j)/2}\sum_{m=0}^{n-j-2l}{n-j-2l
\brack m}_q c_{n-j,\,n-j-2l}(p,q)a^m.
$$

%
%

Based on the orthogonality relation of $q$-Hermite polynomials $H_n(x|q)$
$$
\frac{(q;q)_\infty}{2\pi}\int_0^\pi H_m(x|q)H_n(x|q)(e^{2i\theta},
 e^{-2i\theta};q)_\infty d\theta=(q;q)_n\delta_{mn},
$$
Ismail and Stanton \cite{ismail-stanton03} found two generating
functions for the $q$-Hermite polynomials
\begin{align}
&\sum_{n=0}^\infty\frac{H_{2n}(x|q)}{(q^2;q^2)_n}t^n=\frac{(-t;q)_\infty}{(te^{2i\theta},
 te^{-2i\theta};q^2)_\infty},\label{gfa}\\
&\sum_{n=0}^\infty \frac{H_n(x|q^2)}{(q;q)_n}t^n=\frac{(qt^2;q^2)_\infty}
{(te^{i\theta},te^{-i\theta};q)_\infty}\label{gfb},
\end{align}
where $x=\cos\theta$.

Similarly, from the orthogonality
relation \cite{Koekoek-Swarttouw} of $H_n(x;a|q)$
\begin{equation}\label{ort-big}
\frac{(q;q)_\infty}{2\pi}\int_0^\pi
H_n(x;a|q)H_m(x;a|q)\frac{(e^{2i\theta},
e^{-2i\theta};q)_\infty}{(ae^{i\theta}, ae^{-i\theta};q)_\infty}\
d\theta=(q;q)_n\delta_{mn},
\end{equation}
it follows the identities:
\begin{align}
\sum_{n=0}^\infty \bigg(\sum_{k=0}^{n/2} \frac{q^{n-2k\choose
2}a^{n-2k}t^{n-k}}{(q^2;q^2)_k(q;q)_{n-2k}}\bigg)H_n(x;a|q)&=\frac{(a^2t;q^2)_\infty(-t;q)_\infty}{(te^{2i\theta},
te^{-2i\theta};q^2)_\infty}, \label{gen-1}\\
\sum_{n=0}^\infty
\bigg(\sum_{k=0}^n\frac{(-1)^kq^{k^2}a^kt^{n+k}}{(q^2;q^2)_k(q;q)_{n-k}}\bigg)H_n(x;a|q^2)
&=\frac{(at;q)_\infty(qt^2;q^2)_\infty}{(te^{i\theta},
te^{-i\theta}; q)_\infty},\label{gen-2}
\end{align}
where $x=\cos\theta$.

Clearly,  \eqref{gfa} and \eqref{gfb} are the special cases when $a=0$.
Recall that the  generating function \cite{Koekoek-Swarttouw} of $H_n(x;a|q)$ is
\begin{equation}\label{gfHxa}
\sum_{n=0}^m H_n(x;a|q)\frac{t^n}{(q;q)_n}=\frac{(at;q)_\infty}
{(te^{i\theta},te^{-i\theta};q)_\infty},\quad x=\cos\theta.
\end{equation}

Setting $q\rightarrow p$ and multiplying by the weight function
$(e^{2i\theta}, e^{-2i\theta};q)_\infty/(ae^{i\theta},
ae^{-i\theta};q)_\infty$, the integrals of the generating
functions (4.6--4.8) on base $p$ can be stated as follows
\begin{align*}
& J_{p,q}(a,t)=\frac{(q;q)_\infty(a^2t;p^2)_\infty(-t;p)_\infty}{2\pi}\int_0^\pi
\frac{(e^{2i\theta}, e^{-2i\theta};q)_\infty}{(ae^{i\theta},
ae^{-i\theta};q)_\infty(te^{2i\theta}, te^{-2i\theta};
p^2)_\infty}\ d\theta, \\[7pt]
& H_{p,q}(a,t)=\frac{(q^2;q^2)_\infty(at;p)_\infty(pt^2;p^2)_\infty}{2\pi}\int_0^\pi
\frac{(e^{2i\theta}, e^{-2i\theta};q^2)_\infty}{(ae^{i\theta},
ae^{-i\theta};q^2)_\infty(te^{i\theta}, te^{-i\theta}; p)_\infty}\
d\theta, \\[7pt]
& I_{p,q}(a,t)=\frac{(q;q)_\infty(at;p)_\infty}{2\pi}\int_0^\pi
\frac{(e^{2i\theta}, e^{-2i\theta};q)_\infty}{(ae^{i\theta},
ae^{-i\theta};q)_\infty(te^{i\theta}, te^{-i\theta}; p)_\infty}\
d\theta.
\end{align*}
When $a=0$, they reduce to the integrals $J_{p,q}(t), H_{p,q}(t)$ and $I_{p,q}(t)$
introduced by Ismail and Stanton \cite{ismail-stanton03}.

We conclude  with the observation  that for some special values of $p$ and $q$, the above integrals can be computed
by the Askey-Wilson integral
\cite[p.\,154]{gasperrahman90}:
\begin{align}\label{askey-wilson}
\frac{(q;q)_\infty}{2\pi}\int_0^\pi \frac{(e^{2i\theta},e^{-2i\theta};q)_\infty}
{(ae^{i\theta},ae^{-i\theta},be^{i\theta},be^{-i\theta},ce^{i\theta},ce^{-i\theta},
de^{i\theta},de^{-i\theta};q)_\infty}\ d\theta
=\frac{(abcd;q)_\infty}{(ab,ac,ad,bc,bd,cd;q)_\infty}.
\end{align}

Because of  the change of base formula \eqref{changebase} and
the orthogonality relation \eqref{ort-big},
we may transform the above integrals into summations
\begin{align*}
J_{p,q}(a,t)&=\sum_{n=0}^\infty \sum_{k=0}^{n/2}
\frac{p^{n-2k\choose
2}a^{n-2k}t^{n-k}}{(p^2;p^2)_k(p;p)_{n-2k}}\sum_{j=0}^n d_{n,j,l,n-j-2l}(p,q),\\
H_{p,q}(a,t)&=\sum_{n=0}^\infty
\sum_{k=0}^n\frac{(-1)^kp^{k^2}a^kt^{n+k}}{(p^2;p^2)_k(p;p)_{n-k}}\sum_{j=0}^n d_{n,j,l,n-j-2l}(p^2,q), \\
I_{p,q}(a,t)&=\sum_{n=0}^\infty \frac{t^n}{(p;p)_n} \sum_{j=0}^n
d_{n,j,l,n-j-2l}(p,q).
\end{align*}

Using special cases
of the Askey-Wilson integral \eqref{askey-wilson}, we give some examples of
$H_{p,q}(a,t)$ which have closed product formulas:
\begin{align*}
H_{q,q}(a,t)&=\sum_{n=0}^\infty
\sum_{k=0}^n\frac{(-1)^kq^{k^2}a^kt^{n+k}}{(q^2;q^2)_k(q;q)_{n-k}}\sum_{j=0}^n
d_{n,j,l,n-j-2l}(q^2,q)=1,\\[6pt]
H_{-q,q}(a,t)&=\sum_{n=0}^\infty
\sum_{k=0}^n\frac{(-1)^k(-q)^{k^2}a^kt^{n+k}}{(q^2;q^2)_k(-q;-q)_{n-k}}\sum_{j=0}^n
d_{n,j,l,n-j-2l}(q^2,q)
=1,\\[6pt]
H_{q^2,q}(a,t)&=\sum_{n=0}^\infty
\sum_{k=0}^n\frac{(-1)^kq^{2k^2}a^kt^{n+k}}{(q^4;q^4)_k(q^2;q^2)_{n-k}}\sum_{j=0}^n
d_{n,j,l,n-j-2l}(q^4,q)=(q^2t^2;q^4)_\infty,\\[6pt]
H_{q^2,q^3}(a,t)&=\sum_{n=0}^\infty
\sum_{k=0}^n\frac{(-1)^kq^{2k^2}a^kt^{n+k}}{(q^4;q^4)_k(q^2;q^2)_{n-k}}\sum_{j=0}^n
d_{n,j,l,n-j-2l}(q^4,q^3)=\frac{(at^3q^6;q^6)_\infty}{(t^2q^4;q^4)_\infty}.
\end{align*}

\vspace{.2cm} \noindent{\bf Acknowledgments.} We are grateful to
the referees for valuable suggestions that lead to improvements of
an earlier version.  We would also like to thank Qing-Hu Hou,
Vincent Y. B. Chen and Nancy S. S. Gu for helpful comments. This
work was supported by the 973 Project, the PCSIRT Project of the
Ministry of Education, the Ministry of Science and Technology, and
the National Science Foundation of China.


\end{document}